\documentclass[12pt]{amsart}
\usepackage{amssymb}
\usepackage{amscd}
\usepackage{pstricks,pst-node}
\psset{linewidth=0.4pt}
\input{tableau.tex}

\renewcommand{\H}{\mathcal{H}}
\newcommand{\B}{\mathcal{B}}
\newcommand{\N}{\mathbb{N}}
\newcommand{\cN}{\mathcal{N}}
\newcommand{\F}{\mathcal{F}}
\renewcommand{\O}{\mathcal{O}}

\newcommand{\Z}{\mathbb{Z}}
\newcommand{\C}{\mathbb{C}}
\newcommand{\bd}{\mathbf{d}}
\newcommand{\bl}{{\boldsymbol{\lambda}}}
\newcommand{\bm}{{\boldsymbol{\mu}}}
\newcommand{\bn}{{\boldsymbol{\nu}}}
\newcommand{\br}{{\boldsymbol{\rho}}}
\newcommand{\bL}{\boldsymbol{\Lambda}}
\newcommand{\gl}{\mathfrak{gl}}
\newcommand{\Hom}{\mathrm{Hom}}
\newcommand{\End}{\mathrm{End}}
\newcommand{\Ind}{\mathrm{Ind}}

\numberwithin{equation}{section}

\theoremstyle{plain}
\newtheorem{theorem}{Theorem}[section]
\newtheorem{corollary}[theorem]{Corollary}
\newtheorem{lemma}[theorem]{Lemma}

\theoremstyle{definition}

\newtheorem{example}[theorem]{Example}
\newtheorem{remark}[theorem]{Remark}

\begin{document}
\title{Two-row nilpotent orbits of cyclic quivers}

\author{Anthony Henderson}
\address{School of Mathematics and Statistics,
University of Sydney, NSW 2006, AUSTRALIA}
\email{anthonyh@maths.usyd.edu.au}
\begin{abstract}
We prove that the local intersection cohomology of nilpotent
orbit closures of cyclic quivers is trivial when the two orbits
involved correspond to partitions with at most two rows. This
gives a geometric proof of a result of Graham and Lehrer,
which states that standard modules of the affine Hecke algebra
of $GL_d$ corresponding to nilpotents with at most two Jordan
blocks are multiplicity-free.
\end{abstract}
\maketitle
\section{Introduction}
Let $\Delta_n$ be the cyclic quiver with $n$ vertices. The isomorphism
classes of nilpotent complex representations of $\Delta_n$ are in
bijection with certain \emph{nilpotent orbits} in a suitable variety
(see \S2 for the definition). The closure of such a nilpotent orbit
is usually singular, and it is an important problem to compute its
local intersection cohomology at the points of another given orbit (see \S3).
If $n=1$, these are the usual nilpotent orbits in
$\mathfrak{gl}_d$, and the problem was solved by Lusztig in \cite{greenuni}.

The main result of this paper (Theorem \ref{mainthm} below) is that
this local intersection cohomology is trivial when the two orbits
involved correspond to partitions with at most two rows
(equivalently, the representations of $\Delta_n$ involved are either
indecomposable or the sum of two indecomposables).

In \S3 we discuss some equivalent formulations of this result. The most
noteworthy concerns the complex representation theory of 
$\widetilde{\H}_d$, the affine Hecke algebra of $GL_d(\C)$
(as defined in \cite{kazhdanlusztig}). This algebra 
has geometrically-defined \emph{standard modules} 
$M_{s,x,q_0}$ indexed by triples $(s,x,q_0)$ where $s\in GL_d(\C)$ 
is semisimple,
$x\in\gl_d(\C)$ is nilpotent, $q_0\in\C^{\times}$, and
\[ \mathrm{Ad}(s)(x)=q_0x. \]
By a well-known result of Ginzburg, Theorem \ref{mainthm} (for $n\geq 2$)
is equivalent to the statement that $M_{s,x,q_0}$ is multiplicity-free
when $q_0$ is a primitive $n$-th root of unity, all eigenvalues of
$s$ are powers of $q_0$, and $x$ has at most two Jordan blocks.
(See \S3 for more on the equivalence.)

This statement was recently proved by Graham and Lehrer
in \cite{grahamlehrer}. Actually, they imposed no conditions
on $s$ and $q_0$, but it is easy to reduce to the above case. Their method
is algebraic, making use of the representation theory
of the affine Temperley-Lieb algebras, and their results apply over more
general ground fields (with induced modules instead of standard modules).
The present paper arose from the author's desire to give
a purely geometric/combinatorial proof
of Theorem \ref{mainthm}, thus reproving Graham and Lehrer's result
(in the complex case).

In \S4 we introduce canonical resolutions of the nilpotent orbit
closures, and derive a dimension formula for the orbits, which appears
to be new (see Lemma \ref{geomlemma}). We also discuss how to determine
the local intersection cohomology, given the Poincar\'e polynomials
of the fibres of these resolutions. Then for the remainder of the paper
we restrict attention to the nilpotent orbits with at most two
rows. Section 5 contains some preliminary Lemmas,
and Section 6 completely describes the closure relations among such
orbits. Section 7 presents the proof of Theorem \ref{mainthm}.

Our notation for partitions and Young diagrams follows
\cite{macdonald}. In particular, if $\lambda$ is a partition, its nonzero parts
are $\lambda_1\geq\lambda_2\geq\cdots\geq\lambda_{\ell(\lambda)}$,
where $\ell(\lambda)$ is the length; and $\lambda'$ denotes the transpose
partition.

\textit{Acknowledgements.}
I am very grateful to John Graham and Gus Lehrer for telling me
about their result and showing me a draft of \cite{grahamlehrer}.
\section{Nilpotent Orbits of Cyclic Quivers}
Fix a positive integer $n$. Let $\Delta_n$ be the cyclic
quiver of type $\widetilde{A_{n-1}}$, which has vertex set $I=\Z/n\Z$
and an arrow from $i$ to $i-1$ for all $i\in I$. Let $V=\oplus_{i\in I}V_i$
be an $I$-graded vector space over $\C$. We write $d_i$ for
$\dim V_i$, $\bd$ for the dimension vector
$\textbf{dim}\,V=(d_i)_{i\in I}$, and $d$ for $\dim V=\sum_{i\in I}d_i$.
Let
\[ E_V=\bigoplus_{i\in I}\Hom(V_i,V_{i-1}) \]
be the space of \emph{representations} of $\Delta_n$ on $V$.
Two representations in $E_V$ are isomorphic iff they are in the
same orbit of $G_V=\prod_{i\in I} GL(V_i)$, acting on $E_V$ by conjugation.

We say that $x\in E_V$, and its $G_V$-orbit, are \emph{nilpotent} if
$x$ is nilpotent as an element of $\End(V)$. 
Let $\cN_V$ be the subvariety of nilpotent elements of $E_V$.
The nilpotent $GL(V)$-orbits
in $\End(V)$ are in bijection with $\Lambda(d)$, the set of partitions
of size $d$: let $\O_\lambda$ denote the orbit corresponding to the partition
$\lambda$. The $G_V$-orbits in $\cN_V$ have a similar description, as follows.

The isomorphism classes of indecomposable nilpotent 
representations of $\Delta_n$ are called \emph{segments}. There is
a segment $[i;l)$ for each $i\in I$ and $l\in\Z^{+}$; in our convention
this is the isoclass of indecomposables of length $l$ whose socle is
the simple module corresponding to $i$. The notation is meant to
suggest the multiset $\{i,i+1,\cdots,i+l-1\}$, which (regarded as a
function $I\to\N$) is precisely the dimension vector of an
indecomposable representation in this class. 
Since any representation is a direct
sum of indecomposables, the nilpotent representations are parametrized
by \emph{multisegments} (multisets of segments), which we write with direct
sum notation, e.g.\ $[i_1;l_1)\oplus\cdots\oplus[i_s;l_s)$.

Collecting together all segments $[i;l)$ with the same $i$,
we can think of a multisegment as an $I$-tuple of partitions
$\bl=(\lambda^{(i)})_{i\in I}$. Let $\bL$ be the set of such $I$-tuples,
and for $\bl\in\bL$,
write $\lambda$ for the union $\coprod_{i\in I}\lambda^{(i)}$.
We will identify $\bl$ with the corresponding $I$-labelled Young diagram
of shape $\lambda$. This is the unique $I$-labelled Young diagram
in which the labels from
left to right across a row increase by $1$ at each step,
and $\lambda^{(i)}$ is the subpartition
formed by all rows starting with the label $i$. (Two $I$-labelled
Young diagrams which differ by permuting rows of equal length are
considered the same.)
\begin{example}
Suppose $n=3$. Consider the multisegment
\[ [2;2)\oplus[0;2)\oplus[1;4)\oplus[1;1)\oplus[0;2). \]
The corresponding $\bl$ has $\lambda^{(0)}=(2^{2})$,
$\lambda^{(1)}=(41)$, and $\lambda^{(2)}=(2)$. So
\[ \lambda=(42^{3}1),\text{ and }\bl=
\begin{tableau}
:.1.2.0.1\\
:.0.1\\
:.0.1\\
:.2.0\\
:.1\\
\end{tableau}. \]
\end{example}
The dimension vector of a nilpotent representation in the class corresponding
to $\bl$ is $\bd(\bl)=(d_i(\bl))_{i\in I}$, where $d_i(\bl)$ is 
the number of boxes labelled $i$. 
Consequently, the $G_V$-orbits in $\cN_V$ are in bijection with 
$\bL(\bd)=\{\bl\in\bL\,|\,\bd(\bl)=\bd\}$.
For $\bl\in\bL(\bd)$, write $\O_\bl$ for the corresponding orbit. In other
words,
\[ \O_\bl=\{x\in\cN_V\,|\, \textbf{dim}\, \ker x^{k}=\bd(\bl^{\leq k})\}, \]
where $\bl^{\leq k}$
is obtained from $\bl$ by deleting all but the first $k$ columns.
Clearly the $GL(V)$-orbit containing
$\O_\bl$ is $\O_\lambda$.

We define a partial order $\leq$ on $\bL(\bd)$ by
\[ \bm\leq\bl\Leftrightarrow\O_\bm\subseteq\overline{\O_\bl}. \]
We will be mainly interested in the sub-poset $\bL^{\leq 2}(\bd)$
of $\bL(\bd)$ consisting of all $\bl$ with at most two rows. 
The corresponding nilpotent
representations are either indecomposable or the sum of two
indecomposables; as elements of $\mathrm{End}(V)$, they have at most two
Jordan blocks. Obviously $\bm\in\bL^{\leq 2}(\bd)$, $\bm\leq\bl \Rightarrow
\bl\in\bL^{\leq 2}(\bd)$.

\begin{remark}
All of the above remains true for the linear quiver of type $A_{\infty}$,
if we take $I=\Z$. In this case $\bL^{\leq 2}(\bd)$ has at most $2$ 
elements, of which the maximal one
parametrizes the orbit which is dense in $\cN_V=E_V$. As we will see,
$\bL^{\leq 2}(\bd)$ is more interesting for $\Delta_n$,
but still much simpler to describe than the whole poset $\bL(\bd)$.
\end{remark}
\section{Local Intersection Cohomology}
For $\bl\in\bL(\bd)$, let $\H^{k}IC(\overline{\O_\bl})$ be the
$k$-th intersection cohomology sheaf of the (usually singular)
variety $\overline{\O_\bl}$. As we will see, this vanishes if
$k$ is odd. If $\bm\leq\bl$, let the \emph{local
IC polynomial} $\tilde{K}_{\bl,\bm}(t)\in\N[t]$ be
\[ \tilde{K}_{\bl,\bm}(t)=\sum_{k\geq 0}
\dim\H_{x_\bm}^{2k}IC(\overline{\O_\bl})\, t^{k}, \]
where $\H_{x_\bm}^{2k}IC(\overline{\O_\bl})$ means the stalk at
some point $x_\bm\in\O_\bm$. Trivially we have
$\tilde{K}_{\bl,\bl}(t)=1$ for all $\bl$. 
By basic properties of intersection cohomology, $\tilde{K}_{\bl,\bm}(t)$
has constant term $1$ and satisfies
\begin{equation} \label{degeqn}
\deg \tilde{K}_{\bl,\bm} < (\mathrm{codim}\,_{\overline{\O_\bl}}\,\O_\bm)/2,
\text{ if $\bm<\bl$.}
\end{equation}
If $\bm\nleq\bl$,
we set $\tilde{K}_{\bl,\bm}(t)=0$. The main result of this paper is:
\begin{theorem} \label{mainthm}
For $\bl\geq\bm$ in $\bL^{\leq 2}(\bd)$, $\tilde{K}_{\bl,\bm}(t)=1$.
\end{theorem} \noindent
In other words, for $\bl\geq\bm$ in $\bL^{\leq 2}(\bd)$,
$\overline{\O_\bl}$ is \emph{rationally smooth} at the points of $\O_\bm$
(in general, it is not actually smooth there).

If $n=1$, $\bL(\bd)=\Lambda(d)$ and $\O_\bl=\O_\lambda$.
Lusztig proved in \cite{greenuni} that 
\[ \tilde{K}_{\lambda,\mu}(t)=t^{n(\mu)-n(\lambda)}K_{\lambda,\mu}(t^{-1}), \]
where $K_{\lambda,\mu}(\cdot)$ is the Kostka-Foulkes polynomial
and $n(\lambda)=\sum (k-1)\lambda_k$. So in this case Theorem \ref{mainthm}
is well known.

Several important results involve the polynomials
$\tilde{K}_{\bl,\bm}(t)$, which is why one wants to compute them.
In \cite[\S11]{quivers1} Lusztig showed that
the orbit closures $\overline{\O_\bl}$ could be embedded as open
subvarieties of certain affine Schubert varieties of type $A$,
which explains why $\H^{k}IC(\overline{\O_\bl})=0$ for $k$ odd.
It also means that there is an order-preserving injection 
\[ \bL(\bd)\hookrightarrow \tilde{W}:\bl\mapsto w_\bl, \]
where $\tilde{W}$ is the affine Weyl group of type
$\widetilde{A_{d-1}}$, such that
$\tilde{K}_{\bl,\bm}(t)$
is the affine Kazhdan-Lusztig polynomial $P_{w_\bm,w_\bl}(t)$. 
(This gives a combinatorial algorithm
for computing $\tilde{K}_{\bl,\bm}(t)$ in general, but it is a highly
impractical one.)
So Theorem \ref{mainthm} assserts that
certain very special affine Kazhdan-Lusztig polynomials are trivial. 
One can also interpret $(\tilde{K}_{\bl,\bm}(t))$ as the transition
matrix between two bases of the generic Hall algebra of $\Delta_n$,
an algebra of great importance in the theory of quantum affine 
$\mathfrak{sl}_{n}$ (see for instance \cite[\S3]{ltv}). We
will not say anything further about either of these points of view.

There are also well-known interpretations of (some of) the values
$\tilde{K}_{\bl,\bm}(1)$ in complex representation theory. In the case
$n=1$, $\tilde{K}_{\lambda,\mu}(1)$ is the Kostka number $K_{\lambda,\mu}$,
which gives the multiplicity of the simple representation
$V_{\lambda}$ of the symmetric group $S_d$ in the induced representation
$\Ind_{S_{\mu_1}\times S_{\mu_2}\times\cdots\times S_{\mu_{\ell(\mu)}}}
^{S_d}(\C)$. So in this
case Theorem \ref{mainthm} asserts that an induced representation of the form
$\Ind_{S_{\mu_1}\times S_{\mu_2}}^{S_d}(\C)$
is multiplicity-free, which is again well known.

If on the other hand $n\geq 2$, let $\zeta\in\C^{\times}$ be a primitive
$n$-th root of unity. Let $s\in GL(V)$
be the semisimple element which acts by the scalar $\zeta^{-i}$ on $V_i$. Then
$E_V=\{x\in\End(V)\,|\,\mathrm{Ad}(s)(x)=\zeta x\}$.
Thus for any $x\in \cN_V$, we have an associated standard module
$M_{s,x,\zeta}$ of $\widetilde{\H}_{d}$, the affine Hecke algebra
of $GL_d(\C)$ mentioned in the Introduction. Up to isomorphism,
this depends only on the
$G_V$-orbit of $x$, so we will write it as $M_{\bl}$ for $\bl\in\bL(\bd)$.

Let $L_\bl$ be the quotient of $M_\bl$ defined in 
\cite[Chapter 8]{chrissginzburg}, which is either $0$ or a simple module.
We say that $\bl\in\bL(\bd)$ is \emph{aperiodic} if for all $m$ there exists
some $i\in I$ such that $m$ is not a part of $\lambda^{(i)}$. 
It follows from results of Lusztig that $L_\bl\neq 0$ iff
$\bl$ is aperiodic (see \cite[\S2 and Appendix]{ltv}).
All simple constituents of $M_\bm$ for $\bm\in\bL(\bd)$ are isomorphic
to some such nonzero $L_\bl$, 
and we have the following multiplicity formula due to Ginzburg 
(\cite[Theorem 8.6.15]{chrissginzburg}, \cite[\S2]{ltv}):
\[ [M_\bm:L_\bl]=\tilde{K}_{\bl,\bm}(1). \]
So Theorem \ref{mainthm} implies that if 
$\bm\in\bL^{\leq 2}(\bd)$, $M_\bm$ is multiplicity-free.
In fact, this is equivalent to Theorem \ref{mainthm} (for $n\geq 2$),
since an element of $\bL^{\leq 2}(\bd)$ must be aperiodic, unless $n=2$,
all $d_i$ are equal, and
it is the minimal element of $\bL^{\leq 2}(\bd)$. As mentioned in the
Introduction, this equivalent formulation is part of the main result of 
\cite{grahamlehrer}. Thus Theorem \ref{mainthm}
has already been proved representation-theoretically; we will prove it
geometrically in the course of subsequent sections.
\section{Resolutions of the Orbit Closures}
The key to our approach is the fact that each orbit closure
$\overline{\O_\bl}$ has an obvious resolution of singularities.
Fix $\bl\in\bL(\bd)$. Let $l=\ell(\lambda')$ be the number of columns of $\bl$.
Recall that $\bl^{\leq k}\in\bL$ is obtained by deleting all but
the first $k$ columns of $\bl$. Define the $I$-graded partial flag variety
\[ \F_\bl=\{0=W^{(0)}\subset W^{(1)}\subset\cdots
\subset W^{(l)}=V\,|\, \textbf{dim}\, W^{(k)}=\bd(\bl^{\leq k})\}, \]
and let
\[ \widetilde{\overline{\O_\bl}}=\{(x,(W^{(k)}))\in\cN_V\times\F_\bl\,|\,
x(W^{(k)})\subseteq W^{(k-1)}, 1 \leq k \leq l\}. \]
In the case $n=1$, this is the cotangent bundle of $\F_\bl=\F_\lambda$.
In general we have:
\begin{lemma} \label{geomlemma}
\begin{enumerate}
\item  $\F_\bl$ is a smooth irreducible projective variety, and
\[ \dim \F_\bl=\sum_{i\in I}\left[\binom{d_i}{2}-
\sum_{k=1}^{l}\binom{d_i(\bl^{\leq k})-d_i(\bl^{\leq k-1})}{2}\right]. \]
\item  The second projection $\widetilde{\overline{\O_\bl}}\to\F_\bl$
is a vector bundle, with fibres of dimension
\[ \sum_{i\in I}\left[\binom{d_i}{2}+
\sum_{k=1}^{l}\binom{d_i(\bl^{\leq k})-d_i(\bl^{\leq k-1})+1}{2}\right]
-\epsilon(\bl), \]
where
\begin{equation*}
\begin{split}
\epsilon(\bl)&=\sum_{\text{rows $R$ of $\bl$}} 
d_{\text{last label of $R$}}(\bl^{\leq |R|})\\
&=\sum_{i\in I}\sum_{p=1}^{\ell(\lambda^{(i)})}
d_{i+\lambda_p^{(i)}-1}(\bl^{\leq \lambda_p^{(i)}}).
\end{split}
\end{equation*}
\item  The first projection $p_\bl:\widetilde{\overline{\O_\bl}}\to\cN_V$
has image $\overline{\O_\bl}$ and is an isomorphism over $\O_\bl$.
Hence $p_\bl$ is a resolution of $\overline{\O_\bl}$, and
\begin{equation*}
\begin{split}
\dim \O_\bl &= \dim \overline{\O_\bl} = 
\dim \widetilde{\overline{\O_\bl}}\\
&= \sum_{i\in I}d_i^{2}-\epsilon(\bl).
\end{split}
\end{equation*}
So if $\bm\leq\bl$, $\mathrm{codim}\,_{\overline{\O_\bl}}\,\O_\bm=
\epsilon(\bm)-\epsilon(\bl)$.
\end{enumerate}
\end{lemma}
\begin{proof}
(1) is obvious, because $\F_\bl$ is the product of partial flag varieties
in each $V_i$ of the required dimension. 
Clearly $\widetilde{\overline{\O_\bl}}\to\F_\bl$
is a vector bundle, and the fibre over $(W^{(k)})\in\F_\bl$ is
\[ \bigoplus_{i\in I}\{\phi\in\Hom(V_{i+1},V_i)\,|\,
\phi(W_{i+1}^{(k)})\subseteq W_i^{(k-1)}, 1\leq k \leq l\}, \]
which has dimension
\[ \sum_{i\in I}\sum_{k=1}^{l}(d_{i+1}(\bl^{\leq k})-d_{i+1}(\bl^{\leq k-1}))
d_i(\bl^{\leq k-1}). \]
This is the number of ordered pairs of boxes $(b,b')$ in $\bl$ where
the column containing $b'$ is further right than the column containing $b$,
and the label of $b'$ is one more than the label of $b$. Considering the
box $b''$ immediately to the left of $b'$, we see that this is the same
as the number of ordered pairs of boxes $(b,b'')$ with the same label
where the column containing $b''$ is further right than or equal to
the column containing $b$, and where $b''$ is not at the end of a row.
This gives (2). For (3), $p_\bl$ is proper since $\F_\bl$ is projective.
From the definitions,
\[ p_\bl^{-1}(\O_\bl)=\{(x,(W^{(k)}))\in\cN_V\times\F_\bl\,|\,
\ker x^{k}=W^{(k)}, 1 \leq k \leq l\}, \]
which is clearly open in $\widetilde{\overline{\O_\bl}}$ and isomorphic
to $\O_\bl$. The rest follows since we know from (1) and (2) that
$\widetilde{\overline{\O_\bl}}$ is smooth and irreducible.
\end{proof}

Note that $\sum_{i\in I}d_i^{2}=\dim G_V$. So (3) implies that
for $x_\bl\in\O_\bl$, $\dim Z_{G_V}(x_\bl)=\epsilon(\bl)$. (Thus the notation
agrees with \cite[\S3]{ltv}.) Since $ Z_{G_V}(x_\bl)$ is an open subvariety
of the algebra of endomorphisms of $x_\bl$ as a representation of
$\Delta_n$, it is connected. So the only $G_V$-equivariant
simple perverse sheaves on $\cN_V$ are the shifted intersection
cohomology complexes $IC(\overline{\O_\bl})[\dim \O_\bl]$, for
$\bl\in\bL(\bd)$. By the Equivariant Decomposition Theorem 
(see for example \cite[Theorem 8.4.7]{chrissginzburg}),
\begin{equation} \label{decompeqn}
R(p_\bl)_{*}\C[\dim \O_\bl]\cong
\bigoplus_{\substack{\bn\leq\bl\\j\in\Z}}
IC(\overline{\O_\bn})[\dim \O_\bn + j]^{\oplus a_{\bl,\bn,j}}
\end{equation}
for some $a_{\bl,\bn,j}\in\N$. Since the left-hand side is Verdier self-dual,
$a_{\bl,\bn,j}=a_{\bl,\bn,-j}$. Since $p_\bl$ is an isomorphism
over $\O_\bl$, $a_{\bl,\bl,j}$ is $1$ if $j=0$ and $0$ otherwise.

Taking stalk at $x_\bm\in\O_\bm$ of both sides of \eqref{decompeqn}, we get
\[ \sum_{k\geq 0}\dim H^{k}(p_\bl^{-1}(x_\bm))\, t^{k}
= \sum_{\bl\geq\bn\geq\bm}
\left(\sum_{j\in\Z} a_{\bl,\bn,j}\, t^{\epsilon(\bn)-\epsilon(\bl)-j}\right)
\tilde{K}_{\bn,\bm}(t^{2}). \]
Now
\[ p_\bl^{-1}(x_\bm)\cong\{(W^{(k)})\in\F_\bl\,|\,
x_\bm(W^{(k)})\subseteq W^{(k-1)}, 1 \leq k \leq l\}. \]
It is easy to prove (say by induction on $l$) that
$p_\bl^{-1}(x_\bm)$ has a paving by affine spaces, and hence has no
odd cohomologies. (In the cases we use below we will see this another way.)
So in fact $a_{\bl,\bn,j}=0$ if
$j$ is of opposite parity to $\epsilon(\bn)-\epsilon(\bl)$, and
if we write $g_{\bl,\bm}(t)$ for the Poincar\'e polynomial
$\sum_{k\geq 0}\dim H^{2k}(p_\bl^{-1}(x_\bm))\, t^{k}$, we have
\begin{equation} \label{maineqn}
g_{\bl,\bm}(t)
= \sum_{\bl\geq\bn\geq\bm}
\left(\sum_{j\equiv\epsilon(\bn)-\epsilon(\bl)\text{ mod }2} 
a_{\bl,\bn,j}\, t^{(\epsilon(\bn)-\epsilon(\bl)-j)/2}\right)
\tilde{K}_{\bn,\bm}(t).
\end{equation}
Thanks to the degree constraint and others we have mentioned,
knowing $g_{\bl,\bm}(t)$ for all $\bl\geq\bm$ determines
all $a_{\bl,\bm,j}$ and $\tilde{K}_{\bl,\bm}(t)$. In fact:
\begin{lemma} \label{uniqlemma}
If $b_{\bl,\bm,j}\in\C$ and $L_{\bl,\bm}(t)\in\C[t]$ for all
$\bl\geq\bm\in\bL(\bd)$ and $j\in\Z$ satisfy:
\begin{enumerate}
\item for all $\bl\geq\bm$,
\[ g_{\bl,\bm}(t)
= \sum_{\bl\geq\bn\geq\bm}
\left(\sum_{j\equiv\epsilon(\bn)-\epsilon(\bl)\text{ mod }2} 
b_{\bl,\bn,j}\, t^{(\epsilon(\bn)-\epsilon(\bl)-j)/2}\right)
L_{\bn,\bm}(t), \]
\item  $b_{\bl,\bl,j}=\delta_{j,0}$,
\item  $b_{\bl,\bm,j}=b_{\bl,\bm,-j}$,
\item  $b_{\bl,\bm,j}=0$ if $j$ is of opposite parity to
$\epsilon(\bm)-\epsilon(\bl)$,
\item  $L_{\bl,\bl}(t)=1$, and
\item  $\deg L_{\bl,\bm}<(\epsilon(\bm)-\epsilon(\bl))/2$ for $\bm<\bl$,
\end{enumerate}
then $b_{\bl,\bm,j}=a_{\bl,\bm,j}$, $L_{\bl,\bm}(t)=\tilde{K}_{\bl,\bm}(t)$
for all $\bl\geq\bm$, $j\in\Z$.
\end{lemma}
\begin{proof}
We know that $b_{\bl,\bm,j}=a_{\bl,\bm,j}$, 
$L_{\bl,\bm}(t)=\tilde{K}_{\bl,\bm}(t)$ is a solution of (1)--(6).
We now give an algorithm to determine $b_{\bl,\bm,j}$ and $L_{\bl,\bm}(t)$
from (1)--(6), assuming that all $g_{\bl,\bm}(t)$ are known,
thus showing that this
is the unique solution. This algorithm is by induction on
$\epsilon(\bm)-\epsilon(\bl)$. If this is $0$, $\bm=\bl$,
and $b_{\bl,\bl,j}$ and $L_{\bl,\bl}(t)$ are determined by (2) and (5)
respectively. Otherwise, we can apply the induction hypothesis
to $\bl>\bn$ and $\bn>\bm$, and hence assume that all terms in the
right-hand side of (1) have been determined except $\bn=\bl$ and $\bn=\bm$.
So we know the value of
\[ L_{\bl,\bm}(t)+\sum_{j\equiv\epsilon(\bm)-\epsilon(\bl)\text{ mod }2} 
b_{\bl,\bm,j}\, t^{(\epsilon(\bm)-\epsilon(\bl)-j)/2}. \]
Because of (6), this determines $b_{\bl,\bm,j}$ for $j\leq 0$,
$j \equiv\epsilon(\bm)-\epsilon(\bl)\text{ mod }2$, hence for all $j$
by (3) and (4). Thus $L_{\bl,\bm}(t)$ also is determined.
\end{proof} \noindent
The argument of this proof is familiar in other contexts, for instance
in the study of Kazhdan-Lusztig polynomials.
\begin{remark} \label{intervalrem}
Note that if we want to determine a specific polynomial 
$\tilde{K}_{\bl,\bm}(t)$ by this method, we only need to know
$g_{\br,\bn}(t)$ for those $\br$, $\bn$ such that
$\bl\geq\br\geq\bn\geq\bm$; that is, we can restrict attention to the
interval $[\bm,\bl]\subseteq\bL(\bd)$. In particular, in proving
Theorem \ref{mainthm} we can restrict attention to $\bL^{\leq 2}(\bd)$.
In the remaining sections, we will calculate $g_{\bl,\bm}(t)$ for $\bl\geq\bm$
in $\bL^{\leq 2}(\bd)$ and show that it has the required form.
\end{remark}
\begin{remark}
Note that in the case $n=1$, $p_\lambda:\widetilde{\overline{\O_\lambda}}\to
\overline{\O_\lambda}$ is \emph{semismall}; equivalently,
$a_{\lambda,\nu,j}=0$ if $j\neq 0$. For example, $p_{(d)}$ is the famous
Springer resolution of the nilpotent cone. The Poincar\'e polynomial
$g_{\lambda,\mu}(t)$ is the generalized Hall polynomial
$g_{(1^{\lambda_{l}'}),\cdots,(1^{\lambda_{1}'})}^{\mu}(t)$,
the function
$\epsilon(\lambda)$ is $2n(\lambda)+d$, and the nonzero coefficient
$a_{\lambda,\nu,0}$ is the Kostka number $K_{\nu',\lambda'}$
(see \cite[III.6, Example 5]{macdonald}). 
This can be proved geometrically, using
Fourier transform (see \cite[\S5]{myfourier}). 
When $n\geq 2$, Fourier transform does not seem
to be useful for analysing $R(p_\bl)_{*}\C$.
\end{remark}
\begin{remark}
An example when $p_{\bl}$ is not semismall is the case $n=2$,
$\bl=[0;4)\oplus[0;1)$. If $\bm=[0;2)\oplus[0;2)\oplus[0;1)$, the reader
can check that $a_{\bl,\bm,1}=a_{\bl,\bm,-1}=1$.
\end{remark}
\begin{remark}
When $n=2$, it is easy to see that
\[ \epsilon(\bl)=\frac{(d_0-d_1)^{2}+\epsilon(\lambda)}{2}, \]
whence $\dim\O_\bl=\frac{1}{2}\dim\O_\lambda$. This is well known,
since the orbits $\O_\bl$ are the nilpotent orbits of the symmetric
pair $(\gl_d,\gl_{d_0}\oplus\gl_{d_1})$.
\end{remark}
\begin{remark}
When $n\geq 3$, it is not true in general that
$\epsilon(\bl)$ only depends on $\bd$ and $\lambda$.
(Consider $n=3$, $d_0=2$, $d_1=d_2=1$, $\lambda=(21^{2})$.)
However, it is still true that $\bm\leq\bl$, $\mu=\lambda \Rightarrow
\bm=\bl$. This can be seen by interpreting $\overline{\O_\bl}$ as the image
of $p_\bl$.
\end{remark}
\setlength{\tabwidth}{2.4ex}
\setlength{\tabheight}{2.4ex}
\section{Preliminary Lemmas}
Now we restrict attention to $\bl\in\bL^{\leq 2}(\bd)$ and the corresponding
nilpotent orbits. It is convenient to return to multisegment notation,
and write elements of $\bL^{\leq 2}(\bd)$ as
$[i_1;l_1)\oplus[i_2;l_2)$ where $l_2$ is possibly zero (by convention
$[i;0)$ is the empty partition). For any $m\in\Z$, write $\{m\}$ for the
unique element of $\{0,1,\cdots,n-1\}\cap (m+n\Z)$. So
\[ \lfloor\frac{m}{n}\rfloor=\frac{m-\{m\}}{n},\
\lceil\frac{m}{n}\rceil=\frac{m+\{-m\}}{n}. \]
\begin{lemma} \label{epsilonlemma}
Let $\bl=[i_1;l_1)\oplus[i_2;l_2)$, where the factors are ordered so that
$l_1\geq l_2$. Then
\begin{equation*}
\begin{split}
\epsilon(\bl)&=\lceil\frac{l_1}{n}\rceil+\lceil\frac{l_2}{n}\rceil
+\lceil\frac{l_2-\{i_1-i_2+l_1-1\}}{n}\rceil\\
&\quad\quad\quad+\lceil\frac{l_2-\{i_1-i_2\}}{n}\rceil.
\end{split}
\end{equation*}
\end{lemma}
\begin{proof}
This is trivial from the definition of $\epsilon$.
\end{proof}

Fix $\bl=[i_1;l_1)\oplus[i_2;l_2)$, with $l_1\geq l_2\geq 0$. Let
\begin{equation*}
\begin{split} 
\bL^{\leq 2}(\bd)_{\leq\bl}
&=\{\bm\in\bL^{\leq 2}(\bd)\,|\,\bm\leq\bl\},\\
\bL^{\leq 2}(\bd)_{\leq\bl}^{i_1,i_2}
&=\{\bm=[i_1;m_1)\oplus[i_2;m_2)\,|\,\bm\leq\bl\}.
\end{split}
\end{equation*}
\begin{lemma}
If $l_2\geq 1$, $\bL^{\leq 2}(\bd)_{\leq\bl}=
\bL^{\leq 2}(\bd)_{\leq\bl}^{i_1,i_2}$. If $l_2=0$,
\[ \bL^{\leq 2}(\bd)_{\leq\bl}=\bigcup_{i_2}
\bL^{\leq 2}(\bd)_{\leq\bl}^{i_1,i_2}, \]
and the union is disjoint except for $\bl$.
\end{lemma}
\begin{proof}
If $l_2 \geq 1$, we have $\textbf{dim}\,\ker x=\{i_1,i_2\}$
for any $x\in\O_\bl$, where $\{i_1,i_2\}$ is considered as a
multiset, and hence $\textbf{dim}\,\ker x\geq\{i_1,i_2\}$
for any $x\in\overline{\O_\bl}$. Thus if $\bm\leq\bl$ and 
$\bm\in\bL^{\leq 2}(\bd)$,
$\bm$ must be $[i_1;m_1)\oplus [i_2;m_2)$ for some
$m_1,m_2\in\N$. If $l_2=0$, we know only that for $x\in\overline{\O_\bl}$, 
$\textbf{dim}\,\ker x\geq\{i_1\}$. But if $\bm<\bl$ and 
$\bm\in\bL^{\leq 2}(\bd)$, $\bm$ must have two rows, so
there is a unique $i_2$ such that $\bm=[i_1;m_1)\oplus [i_2;m_2)$.
\end{proof}
\begin{lemma} \label{condlemma}
Suppose $\bm=[i_1;m_1)\oplus [i_2;m_2)\in\bL^{\leq 2}(\bd)
_{\leq\bl}^{i_1,i_2}$. Then 
\[ m_1+m_2=l_1+l_2,\ l_1\geq m_1\geq l_2, \]
and either
\begin{equation*}
\begin{split}
&\textup{(a)}\ m_1\equiv l_1,\ m_2\equiv l_2 \text{ mod }n; \text{ or}\\
&\textup{(b)}\ m_1\equiv i_2-i_1+l_2,\ m_2\equiv i_1-i_2+l_1 \text{ mod }n.
\end{split}
\end{equation*}
(It is possible that \textup{(a)} and \textup{(b)} both hold.)
\end{lemma}
\begin{proof}
That $m_1+m_2=d=l_1+l_2$ is obvious. Since $\O_\bm\subseteq\overline{\O_\bl}$,
we must have $\O_\mu\subseteq\overline{\O_\lambda}$, which implies that
$l_1\geq m_1\geq l_2$. Finally,
\[ \bd([i_1;m_1)\oplus [i_2;m_2))=\bd([i_1;l_1)\oplus [i_2;l_2)). \]
Subtracting $\bd([i_1;m_1)\oplus [i_2;l_2))$ from both sides, we get
\[ \bd([i_2+l_2;m_2-l_2))=\bd([i_1+m_1;l_1-m_1)). \]
This clearly implies either (a) or (b).
\end{proof}

It will be easier to get a converse to this Lemma after
we describe the fibres $p_\bl^{-1}(x_\bm)$. For $s_1,s_2\in\N$, let
$\B_{(s_1,s_2)}$ be the Springer fibre 
$p_{(s_1+s_2)}^{-1}(x_{(s_1)\oplus(s_2)})$,
a special case of this fibre when $n=1$. (The choices of the vector
space of dimension $s_1+s_2$ and the element $x_{(s_1)\oplus(s_2)}\in
\O_{(s_1)\oplus(s_2)}$ are unimportant.)
\begin{lemma} \label{flaglemma}
Let $\bm=[i_1;m_1)\oplus [i_2;m_2)\neq\bl$,
where $m_1$ and $m_2$ satisfy the conditions in Lemma
\ref{condlemma}. Then $p_\bl^{-1}(x_\bm)\cong\B_{(s_1,s_2)}$, where
\[ s_1=\lfloor\frac{m_1-l_2-\{i_2-i_1\}}{n}\rfloor,\
s_2=\lfloor\frac{m_2-l_2}{n}\rfloor, \]
so that
\[ s_1+s_2=s:=\frac{l_1-l_2-\{i_2-i_1\}-\{i_1-i_2+l_1-l_2\}}{n}. \]
(We declare $\B_{(-1,s_2)}$ to be empty.)
\end{lemma}
\begin{proof}
By definition,
\[ p_\bl^{-1}(x_\bm)\cong\{(W^{(k)})\in\F_\bl\,|\,
x_\bm(W^{(k)})\subseteq W^{(k-1)}, 1 \leq k \leq l_1\}. \]
It is clear that for $1\leq k\leq l_2$, $W^{(k)}$ is forced to be
$\ker x_\bm^{k}$. Hence $p_\bl^{-1}(x_\bm)\cong p_{\bl'}^{-1}(x_{\bm'})$ where
\[ \bl'=[i_1+l_2;l_1-l_2),\ \bm'=[i_1+l_2;m_1-l_2)\oplus
[i_2+l_2;m_2-l_2). \]
Since $\bl\neq\bm$, $m_2-l_2>0$. In this new fibre,
if $i_1\neq i_2$, $W^{(1)}$ is forced to be
$(\ker x_{\bm'})_{i_1+l_2}$. Passing to the quotient space and repeating,
we see that $W^{(k)}$ is uniquely determined for $1\leq k\leq\{i_2-i_1\}$, so
$p_{\bl'}^{-1}(x_{\bm'})\cong p_{\bl''}^{-1}(x_{\bm''})$ where
\begin{equation*}
\begin{split}
\bl''&=[i_2+l_2;l_1-l_2-\{i_2-i_1\}),\\ 
\bm''&=[i_2+l_2;m_1-l_2-\{i_2-i_1\})
\oplus [i_2+l_2;m_2-l_2).
\end{split}
\end{equation*}
(If any lengths become negative here we understand that the corresponding
fibre is empty.)

We now wish to apply the dual argument, so note that
the end box of $\bl''$ is labelled $i_1+l_1-1$, and the end labels
of $\bm''$ are $i_1+l_1-1$ and $i_2+l_2-1$ in some order.
If $i_1+l_1\neq i_2+l_2$,
$W^{(l_1-l_2-\{i_2-i_1\}-1)}$ is forced to be the preimage of
$(\mathrm{coker }x_{\bm''})_{i_2+l_2-1}$. Hence
$W^{(l_1-l_2-\{i_2-i_1\}-k)}$ is unique for 
$1\leq k\leq\{l_1-l_2+i_1-i_2\}$,
so $p_{\bl''}^{-1}(x_{\bm''})\cong
p_{\bl'''}^{-1}(x_{\bm'''})$, where
\[ \bl'''=[i_2+l_2;l_1-l_2-\{i_2-i_1\}-\{l_1-l_2+i_1-i_2\})=[i_2+l_2;ns), \]
and $\bm'''$ is obtained from $\bm''$ by deleting the last
$\{l_1-l_2+i_1-i_2\}$ boxes of the row whose end label is
$i_1+l_1-1$. This is the first row if (a) holds and the second row
if (b) holds; if both hold, there are no boxes to be deleted.
In either case,
$\bm'''=[i_2+l_2;ns_1)\oplus [i_2+l_2;ns_2)$.

So in the end it suffices to show that
\[ p_{[i;n(s_1+s_2))}^{-1}(x_{[i;ns_1)\oplus[i;ns_2)})\cong\B_{(s_1,s_2)}. \]
This is obvious, since each $I$-graded flag in the left-hand side
is uniquely determined by the induced flag in the subspace corresponding to 
$i$.
\end{proof}
\begin{example}
To illustrate this proof, here is an example of the sequence of steps
when $n=3$ (the boxes with bold borders are the ones to be removed):
\[ 
\bl=\begin{tableau}
:\b{2}{2}{:.0.1\\:.1.2\\}.2.0.1.2.0.1.2.0\\:;\\
\end{tableau}\, ,\
\bm=\begin{tableau}
:\b{2}{2}{:.0.1\\:.1.2\\}.2.0.1.2.0\\
:;;.0.1.2\\
\end{tableau} \]

\[ \rightsquigarrow
\bl'= \begin{tableau}
:\b{1}{1}{:.2\\}.0.1.2.0.1.2.0\\
\end{tableau}\, ,\
\bm'=\begin{tableau}
:\b{1}{1}{:.2\\}.0.1.2.0\\
:.0.1.2\\
\end{tableau} \]

\[ \rightsquigarrow
\bl''=\begin{tableau}
:.0.1.2.0.1.2\b{1}{1}{:.0\\}\\
\end{tableau}\, ,\
\bm''=\begin{tableau}
:.0.1.2\b{1}{1}{:.0\\}\\
:.0.1.2\\
\end{tableau} \]

\[ \rightsquigarrow
\bl'''=\begin{tableau}
:.0.1.2.0.1.2\\
\end{tableau}\, ,\
\bm'''=\begin{tableau}
:.0.1.2\\
:.0.1.2\\
\end{tableau}\
\rightsquigarrow\B_{(1,1)}. \]
\end{example}
\begin{corollary} \label{condcor}
Let $\bm=[i_1;m_1)\oplus [i_2;m_2)$. Then 
$\bm\in\bL^{\leq 2}(\bd)_{\leq\bl}^{i_1,i_2}$
if and only if 
\[ m_1+m_2=l_1+l_2,\ l_1\geq m_1\geq l_2+\{i_2-i_1\}, \]
and either \textup{(a)} or \textup{(b)} holds.
\end{corollary}
\begin{proof}
We know from Lemma \ref{geomlemma} that $\overline{\O_\bl}$ is the image of
$p_\bl$, so $\bm\leq\bl$ iff $p_\bl^{-1}(x_\bm)$ is nonempty.
Lemma \ref{flaglemma} shows that for $m_1$, $m_2$ as in Lemma \ref{condlemma},
$p_\bl^{-1}(x_\bm)$ is empty if and only if (a) holds and 
$m_1-l_2<\{i_2-i_1\}$.
\end{proof}
\begin{corollary} \label{flagcor}
For $\bm<\bl$,
$p_\bl^{-1}(x_\bm)$ has no odd cohomologies, and
\[ g_{\bl,\bm}(t)
= \sum_{k=0}^{\min\{s_1,s_2\}}(\binom{s}{k}-\binom{s}{k-1})\,t^{k}, \]
where $s_1$, $s_2$ and $s$ are as in Lemma \ref{flaglemma}.
\end{corollary}
\begin{proof}
Springer fibres such as $\B_{(s_1,s_2)}$ have no odd cohomologies.
The right-hand side is the Green polynomial $Q_{(1^{s})}^{(s_1,s_2)}(t)$,
which is well known to equal $\sum_{k\geq 0}
\dim H^{2k}(\B_{(s_1,s_2)})\,t^{k}$.
\end{proof}
\section{Description of the Posets}
Maintain the notations of \S5. In particular,
$\bl=[i_1;l_1)\oplus[i_2;l_2)$ where $l_1\geq l_2\geq 0$; and
$s$ is as in Lemma \ref{flaglemma}.
Using Corollary \ref{condcor},
we can completely describe the poset
$\bL^{\leq 2}(\bd)_{\leq\bl}^{i_1,i_2}$
and hence the poset $\bL^{\leq 2}(\bd)$ as a whole.
The arguments are trivial applications of the Lemmas in \S5, and are
mostly left to the reader. There are four cases.

\textbf{Case 1:} $i_1=i_2=i$, $l_1\equiv l_2\text{ mod }n$.\newline
In this case conditions (a) and (b) coincide, and
$s=\frac{l_1-l_2}{n}$. If $l_1-l_2=0$ or $n$,
$\bl$ is the unique element of $\bL^{\leq 2}(\bd)_{\leq\bl}^{i,i}$, 
and conversely.
Otherwise, $[i;l_1-n)\oplus[i;l_2+n)$ is a predecessor of $\bl$,
also in Case 1. So $\bL^{\leq 2}(\bd)_{\leq\bl}^{i,i}$ is a chain:
\begin{equation*}
\begin{CD}
[i;l_1)\oplus [i;l_2) @. \mathbf{(s,0)}\\
@VVV @.\\
[i;l_1-n)\oplus [i;l_2+n) @. \mathbf{(s-1,1)}\\
@VVV @.\\
: @.\\
@VVV @.\\
\qquad
[i;l_1-\lfloor\frac{s}{2}\rfloor n)\oplus [i;l_2+\lfloor\frac{s}{2}\rfloor n)
\qquad
@. \mathbf{(\lceil\frac{s}{2}\rceil,\lfloor\frac{s}{2}\rfloor)}
\end{CD}
\end{equation*}
Here and in subsequent diagrams $\bm\to\bn$ means that $\bn$ is
a predecessor of $\bm$, and the boldface labels are the values 
of $(s_1,s_2)$ as in Lemma \ref{flaglemma} (whose order
is actually not determined in the case that $i_1=i_2$). 
It is easy to see from
Lemma \ref{epsilonlemma} that
$\epsilon([i;l_1-n)\oplus [i;l_2+n))=\epsilon(\bl)-2$.
So the codimension at each link of the chain is $2$. There are
$\lfloor\frac{s}{2}\rfloor$ links.

\textbf{Case 2:} $i_1=i_2=i$, $l_1\not\equiv l_2\text{ mod }n$.\newline
Here conditions (a) and (b) are mutually exclusive, and 
$s=\lfloor\frac{l_1-l_2}{n}\rfloor$. It is easy to see that
$\bl$ is the unique element of 
$\bL^{\leq 2}(\bd)_{\leq\bl}^{i,i}$ if and only if $l_1-l_2<n$;
otherwise, $[i;l_1-\{l_1-l_2\})\oplus[i;l_2+\{l_1-l_2\})$ 
is a predecessor of $\bl$, also in Case 2, and satisfying condition (b).
Again, $\bL^{\leq 2}(\bd)_{\leq\bl}^{i,i}$ is a chain:
\begin{equation*}
\begin{CD}
[i;l_1)\oplus [i;l_2) @. \mathbf{(s,0)}\\
@VVV @.\\
\qquad[i;l_1-\{l_1-l_2\})\oplus[i;l_2+\{l_1-l_2\})\qquad
@. \mathbf{(s,0)}\\
@VVV @.\\
[i;l_1-n)\oplus [i;l_2+n) @. \mathbf{(s-1,1)}\\
@VVV @.\\
[i;l_1-n-\{l_1-l_2\})\oplus[i;\cdots) @. \mathbf{(s-1,1)}\\
@VVV @.\\
: @.\\
@VVV @.\\
[i;l_1-\lfloor\frac{s}{2}\rfloor n-\{l_1-l_2\})\oplus 
[i;\cdots)
@. \mathbf{(\lceil\frac{s}{2}\rceil,\lfloor\frac{s}{2}\rfloor)}
\end{CD}
\end{equation*}
(To make this and later diagrams legible, some row lengths are omitted; recall
that for every element the sum of the lengths is $l_1+l_2$.)
If $s$ is even, the minimal element is actually the same as
the element which might be thought to be directly above it, namely
$[i;l_1-\frac{s}{2}n)\oplus[i;l_2+\frac{s}{2}n)$. So there are
always $s$ links in the chain. This time
the codimension at each link is $1$.

\textbf{Case 3:} $i_1\neq i_2$, $l_1-l_2\equiv i_2-i_1\text{ mod }n$.\newline
In this case conditions (a) and (b) coincide, and
$s=\lfloor\frac{l_1-l_2}{n}\rfloor$. The condition for
$\bl$ to be the only element of
$\bL^{\leq 2}(\bd)_{\leq\bl}^{i_1,i_2}$ is that $l_1-l_2=\{i_2-i_1\}$. 
Otherwise,
$[i_2;l_1-\{l_1-l_2\})\oplus[i_1;l_2+\{l_1-l_2\})$ is a
predecessor of $\bl$, also in Case 3. The summands here 
are in the right order, i.e.\ $l_1-\{l_1-l_2\}>l_2+\{l_1-l_2\}$.
Once more $\bL^{\leq 2}(\bd)_{\leq\bl}^{i_1,i_2}$ is a chain:
\begin{equation*}
\begin{CD}
[i_1;l_1)\oplus [i_2;l_2) @. \mathbf{(s,0)}\\
@VVV @.\\
\qquad [i_2;l_1-\{l_1-l_2\})\oplus[i_1;l_2+\{l_1-l_2\})
\qquad @. \mathbf{(0,s)}\\
@VVV @.\\
[i_1;l_1-n)\oplus [i_2;l_2+n) @. \mathbf{(s-1,1)}\\
@VVV @.\\
[i_2;l_1-n-\{l_1-l_2\})\oplus[i_1;\cdots) @. \mathbf{(1,s-1)}\\
@VVV @.\\
: @.\\
@VVV @.\\
[i_2;l_1-\lfloor\frac{s}{2}\rfloor n-\{l_1-l_2\})\oplus 
[i_1;\cdots) @. \mathbf{(\lfloor\frac{s}{2}\rfloor,\lceil\frac{s}{2}\rceil)}
\end{CD}
\end{equation*}
If $s$ is even, the minimal element here is also
$[i_1;l_1-\frac{s}{2}n)\oplus[i_2;l_2+\frac{s}{2}n)$, so there
are always $s$ links in the chain. Again the
codimension at each link is $1$.

\textbf{Case 4:} $i_1\neq i_2$, $l_1-l_2\not\equiv i_2-i_1\text{ mod }n$.
\newline
Conditions (a) and (b) are mutually exclusive, and
$s=\lfloor\frac{l_1-l_2}{n}\rfloor$ or $\lfloor\frac{l_1-l_2}{n}\rfloor-1$
depending on whether $\{l_1-l_2\}$ or $\{i_2-i_1\}$ is larger. If
$l_1-l_2<\{i_2-i_1\}$, $\bl$ is the unique element of 
$\bL^{\leq 2}(\bd)_{\leq\bl}^{i_1,i_2}$, and conversely.
Otherwise, write $A=\{i_2-i_1\}$, $B=\{i_1-i_2+l_1-l_2\}$;
then $[i_1;l_1-B)\oplus [i_2;l_2+B)$
and $[i_2;l_1-A)\oplus[i_1;l_2+A)$
are predecessors of $\bl$, also in Case 4, and satisfying condition (b). 
These coincide iff $l_1-l_2<n+\{i_2-i_1\}$, in which
case their common value is a minimal element; otherwise, the
summands are in the right order. The poset
$\bL^{\leq 2}(\bd)_{\leq\bl}^{i_1,i_2}$ is as follows:
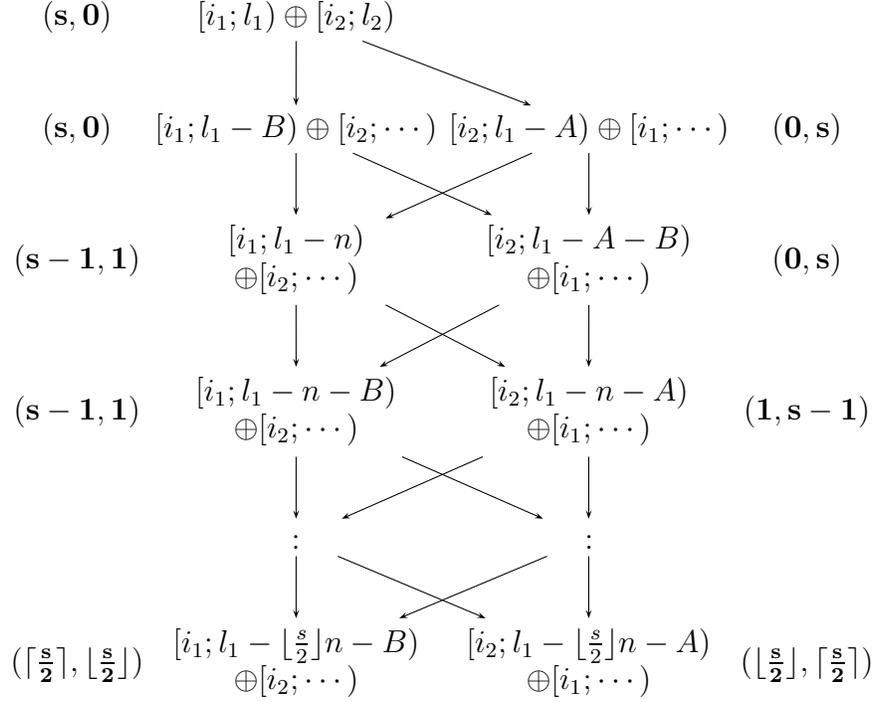
\[
\psmatrix[colsep=0.15cm,rowsep=1cm]
\mathbf{(s,0)} & {[i_1;l_1)\oplus [i_2;l_2)} &  & \\
\mathbf{(s,0)} & {[i_1;l_1-B)\oplus [i_2;\cdots)} & 
{[i_2;l_1-A)\oplus[i_1;\cdots)} & \mathbf{(0,s)}\\
\mathbf{(s-1,1)} & 
{\begin{array}{c}[i_1;l_1-n)\\\oplus[i_2;\cdots)\end{array}} & 
{\begin{array}{c}[i_2;l_1-A-B)\\\oplus[i_1;\cdots)\end{array}} 
& \mathbf{(0,s)}\\
\mathbf{(s-1,1)} & {\begin{array}{c}[i_1;l_1-n-B)\\
\oplus[i_2;\cdots)\end{array}} &
{\begin{array}{c}[i_2;l_1-n-A)\\
\oplus[i_1;\cdots)\end{array}} & \mathbf{(1,s-1)}\\
 &\phantom{AAAAA}:\phantom{AAAAA} & \phantom{AAAAA}:\phantom{AAAAA}& \\
\mathbf{(\lceil\frac{s}{2}\rceil,\lfloor\frac{s}{2}\rfloor)} &
{\begin{array}{c}[i_1;l_1-\lfloor\frac{s}{2}\rfloor n-B)\\
\oplus [i_2;\cdots)\end{array}} &
{\begin{array}{c}[i_2;l_1-\lfloor\frac{s}{2}\rfloor n-A)\\
\oplus [i_1;\cdots)\end{array}} &
\mathbf{(\lfloor\frac{s}{2}\rfloor,\lceil\frac{s}{2}\rceil)}
\endpsmatrix
\psset{nodesep=3pt,arrows=->}
\ncline{1,2}{2,2}
\ncline{1,2}{2,3}
\ncline{2,2}{3,2}
\ncline{2,2}{3,3}
\ncline{2,3}{3,2}
\ncline{2,3}{3,3}
\ncline{3,2}{4,2}
\ncline{3,2}{4,3}
\ncline{3,3}{4,2}
\ncline{3,3}{4,3}
\ncline{4,2}{5,2}
\ncline{4,2}{5,3}
\ncline{4,3}{5,2}
\ncline{4,3}{5,3}
\ncline{5,2}{6,2}
\ncline{5,2}{6,3}
\ncline{5,3}{6,2}
\ncline{5,3}{6,3}
\]
If $s$ is even, the bottom two elements coincide and are minimal.
If $s$ is odd, there is one more element below these two, namely
\[ [i_1;l_1-\lceil\frac{s}{2}\rceil n)\oplus
[i_2;l_2+\lceil\frac{s}{2}\rceil n)\qquad
\mathbf{(\lfloor\frac{s}{2}\rfloor,\lceil\frac{s}{2}\rceil)}. \]
So the length of a maximal chain in this poset is $s+1$. Again, all
predecessors have codimension $1$.

\begin{example}
As an example of how these posets
$\bL^{\leq 2}(\bd)_{\leq\bl}^{i_1,i_2}$ fit together to
make $\bL^{\leq 2}(\bd)$, the following is the case
$n=2$, $d_0=d_1=3$:
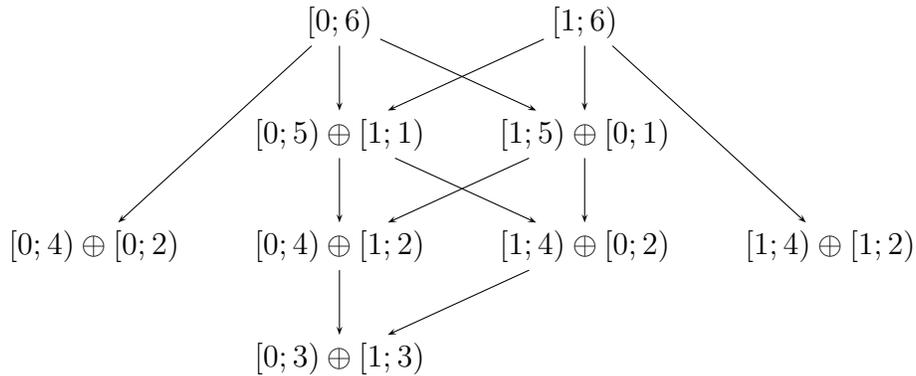
\[
\psmatrix[colsep=1cm,rowsep=1cm]
 & {[0;6)} & {[1;6)} & \\
 & {[0;5)\oplus[1;1)} & {[1;5)\oplus[0;1)} & \\
{[0;4)\oplus[0;2)} & {[0;4)\oplus[1;2)} & {[1;4)\oplus[0;2)} & 
{[1;4)\oplus[1;2)} \\
 & {[0;3)\oplus[1;3)} & &
\endpsmatrix
\psset{nodesep=3pt,arrows=->}
\ncline{1,2}{2,2}
\ncline{1,2}{2,3}
\ncline{1,2}{3,1}
\ncline{1,3}{2,2}
\ncline{1,3}{2,3}
\ncline{1,3}{3,4}
\ncline{2,2}{3,2}
\ncline{2,2}{3,3}
\ncline{2,3}{3,2}
\ncline{2,3}{3,3}
\ncline{3,2}{4,2}
\ncline{3,3}{4,2}
\]
\end{example}
\section{Proof of Theorem \ref{mainthm}}
For $\bl\geq\bm$ in $\bL^{\leq 2}(\bd)$,
define a nonnegative integer $c_{\bl,\bm}$ by the following rule.
If $\bl=\bm$, let $c_{\bl,\bm}=1$. If $\bl>\bm$, then
$\bm\in\bL^{\leq 2}(\bd)_{\leq\bl}^{i_1,i_2}$ for unique $(i_1,i_2)$,
and the interval $[\bm,\bl]$ is contained in 
$\bL^{\leq 2}(\bd)_{\leq\bl}^{i_1,i_2}$.
Define $s_1$, $s_2$, and $s$ as above, and let $k(\bm)=\min\{s_1,s_2\}$.
Say $\bm$ is \emph{special} (relative to $\bl$) if there is no
element $\bn$ of $[\bm,\bl]$ with $k(\bn)=k(\bm)$ except $\bm$ itself.
If $\bm$ is special, let 
$c_{\bl,\bm}=\binom{s}{k(\bm)}-\binom{s}{k(\bm)-1}$.
Otherwise, let $c_{\bl,\bm}=0$. 
We can now prove a stronger
version of Theorem \ref{mainthm}:
\begin{theorem}
If $\bl\geq\bm$ in $\bL^{\leq 2}(\bd)$, then $\tilde{K}_{\bl,\bm}(t)=1$,
$a_{\bl,\bm,j}=0$ if $j\neq 0$, and $a_{\bl,\bm,0}=c_{\bl,\bm}$.
\end{theorem}
\begin{proof}
By inspection of the diagrams in the previous section,
we see that if $\bm$ is special, then $\epsilon(\bm)-\epsilon(\bl)=2k(\bm)$.
Moreover, for all positive $k'\leq k(\bm)$, there is a unique
special $\bn\in[\bm,\bl]$ such that $k(\bn)=k'$. So
Corollary \ref{flagcor} can be rephrased:
\begin{equation*}
g_{\bl,\bm}(t)
= \sum_{\bl\geq\bn\geq\bm}c_{\bl,\bn}\,t^{(\epsilon(\bn)-\epsilon(\bl))/2}.
\end{equation*}
Now we apply Lemma \ref{uniqlemma} (bearing in mind Remark \ref{intervalrem})
with $b_{\bl,\bm,j}=\delta_{j,0}c_{\bl,\bm}$ and $L_{\bl,\bm}(t)=1$. 
The result follows.
\end{proof}

\end{document}